\documentclass[11pt]{amsart}

\usepackage{amsmath,amsthm,amssymb,amscd}

\input{xy}\xyoption{all}
\CompileMatrices

\newcommand{\p}{{\mathfrak{p}}}

\newcommand{\Q}{{\mathbb{Q}}}
\newcommand{\m}{{\mathfrak{m}}}

\newcommand{\Hom}{\operatorname{Hom{}}}
\newcommand{\Ext}{\operatorname{Ext{}}}

\newcommand{\Ann}{\operatorname{Ann}}
\newcommand{\rank}{\operatorname{rank{}}}
\newcommand{\grade}{\operatorname{grade}}

\newcommand{\syz}{\operatorname{syz{}}}
\renewcommand{\hat}{\widehat}
\renewcommand{\bar}{\overline}

\renewcommand{\phi}{\varphi}
\renewcommand{\to}{{\longrightarrow}}
\newcommand{\xx}{{\underline{x}}}

\newcommand{\tr}{\operatorname{tr}}
\newcommand{\ev}{\operatorname{ev}}
\newcommand{\im}{\operatorname{im}}
\newcommand{\coker}{\operatorname{coker}}
\newcommand{\pd}{\operatorname{pd}}

\newtheorem{thm}{Theorem}
\newtheorem{cor}[thm]{Corollary}
\newtheorem*{cor*}{Corollary}
\newtheorem{prop}[thm]{Proposition}
\newtheorem{lemma}[thm]{Lemma}
\newtheorem{defn}[thm]{Definition}

\newtheorem*{conj*}{Conjecture}

\newtheorem*{main}{Main Theorem}

\begin{document}

\title{On a conjecture of Auslander and Reiten}

\author{Craig Huneke}
\address{Department of Mathematics \\
        University of Kansas \\
        Lawrence, KS
        66045}
\email{huneke@math.ukans.edu}
\urladdr{http://www.math.ukans.edu/\textasciitilde huneke}

\author{Graham J. Leuschke}
\curraddr{MSRI\\
1000 Centennial Drive, \#5070\\
       Berkeley, CA 94720-5070}
\email{gleuschke@math.ukans.edu}
\urladdr{http://www.leuschke.org/}

\date{\today}

\thanks{The first author was supported by NSF grant DMS-0098654, and the second author was supported by an NSF Postdoctoral Fellowship and the Clay Mathematics Institute.}

\bibliographystyle{amsplain}

\numberwithin{thm}{section}

\begin{abstract}  In studying Nakayama's 1958 conjecture on rings of
infinite dominant dimension, Auslander and Reiten proposed the following
generalization: Let $\Lambda$ be an Artin algebra and $M$ a
$\Lambda$-generator such that $\Ext^i_\Lambda(M,M)=0$ for all $i \geq
1$; then $M$ is projective.  This conjecture makes sense for any ring. 
We establish Auslander and Reiten's conjecture for excellent
Cohen--Macaulay normal domains containing the rational numbers, and
slightly more generally.\end{abstract}

\maketitle

\setcounter{section}{-1}

\section{Introduction} 

The generalized Nakayama conjecture of M.~Auslander and I.~Reiten is as
follows \cite{AusReit:1975}: For an Artin algebra $\Lambda$,  every
indecomposable injective $\Lambda$-module appears as a direct summand in
the  minimal injective resolution of $\Lambda$. Equivalently, if $M$ is
a finitely generated $\Lambda$-generator such that $\Ext^i_\Lambda(M,M) 
=0$ for all $i\geq 1$, then $M$ is projective.  This latter formulation
makes sense for any ring, and Auslander, S. Ding, and \O. Solberg
\cite{AusDingSolberg} widened the context to algebras over commutative
local rings. 

\begin{conj*}[AR]\label{mainconj} Let $\Lambda$ be a Noetherian ring finite over its center
and $M$ a finitely generated left $\Lambda$-module such that
$\Ext_\Lambda^i(M,\Lambda) = \Ext_\Lambda^i(M,M) = 0$ for all  $i>0$. 
Then $M$ is projective.\end{conj*}

In the same paper, Auslander and Reiten proved AR for modules $M$ that
are {\it ultimately closed}, that is, there is some syzygy $N$ of $M$
all of whose indecomposable direct summands already appear in some
previous syzygy of $M$.  This includes all modules over rings of finite
representation type, all rings $\Lambda$ such that for some integer $n$,
$\Lambda$ has only a finite number of indecomposable summands of
$n^{\text{th}}$ syzygies, and all rings of radical square zero. 

Auslander, Ding, and Solberg \cite[Proposition 1.9]{AusDingSolberg}
established AR in case $\Lambda$ is a quotient of a ring $\Gamma$ of
finite global dimension by a regular sequence.  In fact, in this case
they prove something much stronger: If $\Ext_{\Lambda}^2(M,M)=0$, then
$\pd_\Lambda M < \infty$ \cite[Proposition 1.8]{AusDingSolberg}.  This in turn
was generalized by L. Avramov and R.-O. Buchweitz \cite[Theorem
4.2]{Avramov-Buchweitz}: A finite module $M$ over a (commutative)
complete intersection ring $R$ has finite projective dimension if and
only if $\Ext_R^{2i}(M,M)=0$ for some $i>0$.

M. Hoshino \cite{Hoshino} proved that if $R$ is a symmetric Artin algebra with radical cube zero, then
$\Ext_R^1(M,M)=0$ implies that $M$ is free.
Huneke, L.M. \c Sega, and A.N. Vraciu have recently extended this to prove
that if $R$ is Gorenstein local with $\m^3 = 0$, and if
$\Ext_R^i(M,M)=0$ for some $i \geq 1$, then $M$ is free, and have
further verified the Auslander-Reiten conjecture for
all finitely generated modules $M$ over Artinian commutative local rings $(R,
\m)$ such that
$\m^2M = 0$ \cite{Huneke-Sega-Vraciu}. In  particular, this verifies the
Auslander-Reiten conjecture for commutative local rings with $\m^3 = 0$.

The assumption that $\Lambda$ be finite over its center is essential, given a counterexample due to R. Schultz \cite{Schultz}.

Our main theorem establishes the AR conjecture for a class of commutative
Cohen--Macaulay rings and well-behaved modules.  Moreover, our result is
effective; we can specify
how many $\Ext$ are needed to vanish to give the conclusion of AR.

\begin{main}\label{maintheorem}  Let $R$ be a Cohen--Macaulay ring  which is a quotient of a locally excellent ring $S$ of dimension $d$ by a
locally regular sequence.
Assume that $S$ is 
locally a complete intersection ring in codimension one, 
and further assume either that $S$ is Gorenstein, or that $S$
contains the field of rational numbers.  
Let $M$ be a finitely generated $R$-module of constant rank such that
\begin{equation}\label{ext}\begin{split}
Ext_R^i(M,M)&=0 \quad\text{for $i=1, \dots, d,$ and} \\ 
Ext_R^i(M,R)&= 0 \quad\text{for $i= 1, \dots, 2d+1.$}
\end{split}\end{equation}
Then $M$ is projective.\end{main}

The restriction imposed on $R$ by assuming that $S$ be locally complete intersection in codimension one is equivalent to assuming that $R$ is a quotient by a regular sequence of some normal domain $T$, by \cite[Theorem 3.1]{Huneke-Ulrich:1985}.  However, replacing $S$ by $T$ according to the construction in \cite{Huneke-Ulrich:1985} would increase $d$, the number of $\Ext$ required to vanish.  In any case, this observation gives the following corollary.

\begin{thm}  Let $R$ be a Cohen--Macaulay ring  which is a quotient of a locally excellent ring $S$ of dimension $d$ by a
locally regular sequence.
Assume that $S$ is 
locally a complete intersection ring in codimension one, 
and further assume either that $S$ is Gorenstein, or that $S$
contains the field of rational numbers.    Then the AR conjecture holds for all finitely generated $R$-modules, that is, if $\Ext_R^i(M,R) = \Ext_R^i(M,M) = 0$ for all  $i>0$, then $M$ is projective.\end{thm}

Not every zero-dimensional ring $R$ is a factor of a ring $S$ as in 
the theorem, since not all Artinian local rings can
be smoothed.  For example, Anthony Iarrobino has pointed out that the easiest such example is a polynomial
ring in four variables modulo an ideal generated by seven general quadrics (note, however, that the cube of the maximal ideal of such a ring is zero, so this case is covered by \cite{Huneke-Sega-Vraciu}). For other examples of non-smoothable rings, see Mumford \cite{Mumford}.

In the next section we prove some preliminary lemmas, and then prove the main
result. This requires extra work regarding the trace of a module. Since
we could not find a satisfactory reference for what we needed,
we include basic facts concerning the trace in an appendix.  

Throughout the following, all rings are Noetherian and all modules
finitely generated.  For an $R$-module $M$, we define the {\it dual} of
$M$ by $M^* = \Hom_R(M,R)$.  There is a  natural homomorphism $\theta_M
: M \to M^{**}$ defined by sending $x\in M$ to ``evaluation at $x$''. 
We say that $M$ is {\it torsion-free} if $\theta_M$ is injective, and
{\it reflexive} if $\theta_M$ is an isomorphism.  It is known (cf.
\cite[Theorem 2.17]{AusBridger}, for example) that $M$ is torsion-free
if and only if $M$ is a first syzygy, and reflexive if and only if
 $M$ is a second syzygy.  We
will say that a torsion-free $R$-module $M$ {\it has constant rank} if
$M$ is locally free of constant rank at the minimal primes of $R$.  This
is equivalent to $K \otimes_R M$ being a free $K$-module, where $K$ is
the total quotient ring of $R$ obtained by inverting all
nonzerodivisors.

\section{Proof of the Main Theorem}\label{genreductions}

We begin by observing that the vanishing of $\Ext$ and the projectivity of
$M$ are both local questions, so that in proving our main theorem we may
assume that both $S$ and $R$ are local.  Furthermore, since $S$ is
assumed to be excellent we can (and do) complete $S$ at its maximal
ideal without loss of generality.

Next we point out the following consequence of the lifting criterion of
Auslander, Ding, and Solberg \cite[Proposition 1.6]{AusDingSolberg}. 

\begin{lemma}\label{lift} Let $S$ be a complete local ring, $x \in S$ a
nonunit nonzerodivisor, and $R=S/(x)$.  Assume that there exists $t \geq 2$ such that for any $S$-module $N$, $\Ext_S^i(N,N)=\Ext_S^i(N,S)=0$ for $i=1, \ldots, t$ implies that $N$ is free.  Then for any $R$-module $M$, $\Ext_R^i(M,M)=\Ext_R^i(M,R)=0$ for $i=1, \ldots, t$ implies that $M$ is free.  Furthermore, if AR holds for $S$-modules then it holds for $R$-modules. \end{lemma}

\begin{proof} Let $M$ be an $R$-module such that $\Ext_R^i(M,M)=\Ext_R^i(M,R)=0$ for $i=1, \ldots, t$.  Then in particular $\Ext_R^2(M,M)=0$, and so by
\cite[Proposition 1.6]{AusDingSolberg} there exists an $S$-module $N$ on which
$x$ is a nonzerodivisor and such that $R \otimes_S N \cong M$.  Apply
$\Hom_S(-, N)$ to the short exact sequence $0 \to N \to N \to M \to 0$
and use the fact that $\Ext_S^{i+1}(M,N) \cong \Ext_R^i(M,M) =0$ for $i
=1, \ldots, t$ to see that multiplication by $x$ is surjective on
$\Ext_S^i(N,N)$ for $i=1, \ldots, t$.  Then Nakayama's Lemma implies that
$\Ext_S^i(N,N)=0$ for $i=1, \ldots, t$.  The same argument, applying
$\Hom_S(-,S)$ and observing that $\Ext_S^{i+1}(M,S) \cong \Ext_R^i(M,R)=0$ 
for $i =1, \ldots, t$, 
shows that $\Ext_S^i(N,S)=0$ for $i=1,\ldots, t$ as well.  Since this forces $N$ to be $S$-free, $M$ is $R$-free.

Finally, repeating the argument with ``all $i\geq 1$'' in place of ``$i=1, \ldots t$'' gives the last statement.\end{proof}

With Lemma~\ref{lift} in mind, we now focus on the case $R=S$ in our
main theorem.  Indeed, if $\dim(S) \leq 1$, then $S$ is locally a complete intersection ring by hypothesis, and hence $R$ is as well.  By \cite[Proposition~1.9]{AusDingSolberg}, then, AR holds for $R$-modules.  So we may assume that $R=S$, and in particular we take $d=\dim R$.  
Our next goal is to modify the module $M$.

\begin{lemma}\cite[Lemma 1.4]{AusReit:1975}\label{syzygy}  In proving the Main Theorem, we may replace $M$ by $\syz_R^n(M)$, where $n = \max\{2,
d+1\}$, and assume that $M$ is reflexive and that $\Ext_R^i(M^*,
R)=0$ for $i = 1, \dots, d$.   In proving
AR, we may replace $M$ by any syzygy module $\syz_R^t(M)$.  
\end{lemma}

\begin{proof} Put $N = \syz_R^n(M)$.  It is a straightforward
computation with the long exact sequences of $\Ext$ to show that if
$\Ext_R^i(M,M)=0$ for $i=1, \ldots, d$ and $\Ext_R^i(M,R)=0$ for $i=1, \ldots, 2d+1$, then
$\Ext_R^i(N,N)=0$ for $i=1, \ldots, d$ and $\Ext_R^i(N,R)=0$ for $i=1, \ldots, d$.
Assume, then, that we have shown that $N$ is free.  Then since
$\Ext_R^n(M,R)=0$, the $n$-fold extension of $M$ by $N$
consisting of the free modules in the resolution of $M$ must split, so
$M$ is free as well.  This proves the last statement.

To prove that $N$ is reflexive and $\Ext_R^i(N^*,
R)=0$ for $i = 1, \dots, d$, one shows by induction on $t$ that
$\Ext_R^i((\syz_R^t(M))^*, R)=0$ for $i=1, \dots, t$.  For the base case 
$t=2$, observe that since $\Ext_R^i(M, R)=0$ for $i=1, 2$, the dual of 
the exact sequence
\begin{equation}\label{secondsyz}0 \to \syz_R^2(M) \to F_1 \to F_0 \to M \to 0, \tag{*}\end{equation}
where $F_1$ and $F_0$ are free modules, is still exact.  Dualizing again gives $(\ref{secondsyz})$ back, so $N=\syz_R^2(M)$ is reflexive and satisfies $\Ext_R^i(N^*,R)=0$ for $i=1, 2$.  For the inductive step, dimension-shifting shows that if $\Ext_R^i(M^*,R)=0$ for $i=1, \dots, t-1$, then $\Ext_R^i((\syz_R^1(M))^*,R)=0$ for $i=2, \dots, t$, and the same argument as above shows that $\Ext_R^1((\syz_R^1(M))^*,R)=0$.
\end{proof}

It is worth noting that if $R$ is a Cohen--Macaulay (CM) ring, then
$\syz_R^d(M)$ is a maximal Cohen--Macaulay (MCM) module for any $M$. 
Also, the replacement in Lemma~\ref{syzygy} has consequences for the
assumptions (\ref{ext}) in the main theorem: If $\Ext_R^i(M,R)=0$ for $i=1, \dots, t$, then
$\Ext_R^i(\syz_R^1(M),R)=0$ for $i=1, \dots, t-1$.  This observation
combines with Lemmas~\ref{lift} and \ref{syzygy} to reduce the proof of
our main theorem to the following:

\begin{thm}\label{reducedmain} Let $(R, \m)$ be a complete local CM ring
of dimension $d$ which is a complete intersection in codimension one.
Assume either that $R$ is Gorenstein, or
that $R$ contains $\Q$.  Let $M$ be a MCM $R$-module of constant rank 
such that for $i = 1, \dots, d$,
\begin{equation}\label{exts-d}
\begin{split}
&\Ext_R^i(M,M) = 0,\\
&\Ext_R^i(M,R) = 0, \quad\text{and}\\
&\Ext_R^i(M^*,R) = 0.
\end{split}\end{equation}
Then $M$ is free.
\end{thm}

We postpone the proof of Theorem~\ref{reducedmain} to the end of this section, and establish some preparatory results.

By Cohen's structure theorem, the complete local ring $R$ is a homomorphic image of a regular local ring, and so has a canonical 
module $\omega$.  Since $R$ 
is complete intersection in codimension one, it is in particular 
Gorenstein at the associated primes, and so $\omega$ has constant 
rank. Hence $\omega$ is isomorphic to an ideal of $R$. 
For a MCM $R$-module $N$, we write $N^\vee$ for the canonical dual
$\Hom_R(N,\omega)$.

We next apply a result found in 
\cite[Corollary B4]{Avramov-Buchweitz-Sega} (see also
\cite[Lemma 2.1]{Hanes-Huneke:2002}):

\begin{prop}\label{ABSprop}  Let $R$ be a CM local ring with a canonical
module $\omega$ and let $N$ be a MCM $R$-module.  If $\Ext_R^i(N,R)=0$
for $i=1,\dots,\dim R$ then $\omega\otimes_R N \cong (N^*)^\vee$ is a
MCM $R$-module.\end{prop}

Applied to our current context, this gives the following fact.

\begin{cor} Under our assumptions (\ref{exts-d}) in Theorem~\ref{reducedmain},
both $\omega \otimes_R
M$ and $\omega\otimes_R M^*$ are MCM $R$-modules.\end{cor}

We will also show that the triple tensor product $\omega \otimes_R M^*
\otimes_R M$ is MCM, but for this we use the following lemma.  It
requires that we add one further assumption to (\ref{exts-d}): that the
module $M$ in question has constant rank.  

\begin{lemma}\label{cutdown} Let $(R, \m, k)$ be a CM local ring with
canonical ideal $\omega$, and let $N$ be a MCM $R$-module of constant
rank.  Assume that $\Hom_R(N,N)$ is also a MCM $R$-module, and that for
some maximal regular sequence $\xx$, we have 
$$\Hom_R(N,N) \otimes_R R/(\xx) \cong \Hom_{R/(\xx)}(N/\xx N, N/\xx
N).$$ 
Then $\xx$ is a regular sequence on $N\otimes_R N^\vee$.  In particular,
$N\otimes_R N^\vee$ is MCM.\end{lemma}

\begin{proof}  We indicate reduction modulo $\xx$ by an overline, and
use $\lambda(-)$ for the length of a module.  We also continue to use
$-^\vee$ for $\Hom_{\bar{R}}(-, \bar\omega)$ without fear of confusion. 
Since $\bar\omega \cong E_{\bar{R}}(k)$, the injective hull of the
residue field of $\bar{R}$, we have $\lambda(M^\vee) = \lambda(M)$ for
all $\bar{R}$-modules $M$.

First, a short computation using Hom-Tensor adjointness:
\begin{equation*}\begin{split}
(\bar{N} \otimes_{\bar{R}} \bar{N}^\vee)^\vee &= \Hom(\bar{N}
\otimes_{\bar{R}} \bar{N}^\vee, \bar\omega) \\
        &\cong \Hom_{\bar{R}}(\bar{N}, \bar{N}^{\vee \vee}) \\
        &\cong \Hom_{\bar{R}}(\bar{N}, \bar{N})
\end{split}\end{equation*}
In particular, this implies that $\lambda(\bar{N} \otimes_{\bar{R}}
\bar{N}^\vee) = \lambda((\bar{N} \otimes_{\bar{R}} \bar{N}^\vee)^\vee =
\lambda(\Hom_{\bar{R}}(\bar{N}, \bar{N}))$.  Since $\bar{N\otimes_R N^\vee} = \bar{N} \otimes_{\bar{R}} \bar{N}^\vee$, our hypothesis yields $\lambda(\bar{N \otimes_R N^\vee}) =
\lambda(\bar{\Hom_R(N,N)})$.  Finally, we compute, using the fact that
$N$, $N\otimes_R N^\vee$, and $\Hom_R(N,N)$ all have constant rank:
\begin{equation*}\begin{split}
\lambda(\bar{N\otimes_R N^\vee}) &= \lambda(\bar{\Hom_R(N,N)}) \\
        &= e(\xx, \Hom_R(N,N)) \\
        &= \rank(\Hom_R(N,N)) e(\xx, R) \\
        &= \rank(N)^2 e(\xx, R) \\
        &= \rank(N\otimes_R N^\vee) e(\xx, R) \\
        &= e(\xx, N\otimes_R N^\vee)
\end{split}\end{equation*}
Here $e(\xx, \ )$ denotes the multiplicity of the ideal $(\xx)$ on the
module. The second equality follows since we have assumed that
$\Hom_R(N,N)$ is also a MCM $R$-module.
The equality of the first and last items implies that $N\otimes_R N^\vee$
is MCM by \cite[4.6.11]{BH}.\end{proof}

\begin{prop}\label{simpletriple} Let $(R, \m)$ be a CM local ring with
canonical ideal $\omega$ and let $M$ be a reflexive $R$-module of
constant rank such that $\Ext_R^i(M,M) = \Ext_R^i(M^*,R)=0$ for $i=1,
\dots, d=\dim R$.   Then $\omega \otimes_R M^* \otimes_R M$ is a MCM
$R$-module.\end{prop}

  \begin{proof} We will take $N=M$ in Lemma~\ref{cutdown}.  By
Proposition~\ref{ABSprop}, $M^\vee \cong \omega \otimes_R M^*$, so we need
only show that $\Hom_R(M,M)$ cuts down correctly.   Induction on
the length of a regular sequence $\xx$, using the vanishing of
$\Ext_R^i(M,M)$, then proves that $\xx$ is also regular on $\Hom_R(M,M)$ and that $\Hom_R(M,M) \otimes_R R/(\xx) \cong
\Hom_{R/(\xx)}(M/\xx M, M/\xx M),$ finishing the proof.
 \end{proof}

  \begin{prop}\label{triple}  In addition to the assumptions
(\ref{exts-d}) of Theorem~\ref{reducedmain}, suppose also that $M$ has constant rank. 
Then $\omega \otimes_R M^* \otimes_R M$ is a MCM $R$-module. 
Furthermore, the natural homomorphism 
$$1 \otimes \alpha : \omega \otimes_R M^* \otimes_R M \to \omega
\otimes_R \Hom_R(M,M),$$
where $\alpha$ is defined by $\alpha(f \otimes x)(y) = f(y) \cdot x$, is
injective.\end{prop}

\begin{proof} The first statement follows immediately from
Proposition~\ref{simpletriple}.  For the second, pass to the total quotient
ring $K$ of $R$.  Since $R$ is generically Gorenstein, $\omega \otimes_R K
\cong K$, and since $M$ has a rank, $M \otimes_R K$ is a free $K$-module. 
Since $\alpha$ is an isomorphism when $M$ is free, the kernel of $1
\otimes \alpha$ must be torsion.  But  $\omega \otimes_R M^*\otimes_R M$ is
MCM, and so torsion-free. Hence the kernel of $1 \otimes \alpha$ is
zero.\end{proof}

We return to the assumptions of Theorem~\ref{reducedmain}: $(R, \m, k)$
is a complete local CM ring with a canonical ideal $\omega$, and $M$ is
a torsion-free $R$-module  of constant rank,
 satisfying 
\begin{equation}\label{exts-d-2}
\begin{split}
&\Ext_R^i(M,M) = 0,\\
&\Ext_R^i(M,R) = 0, \quad\text{and}\\
&\Ext_R^i(M^*,R) = 0, \quad\text{for $i = 1, \dots, d=\dim R$.}
\end{split}\end{equation}
We also assume that $R$ is locally a complete intersection ring in
codimension one.  As we observed above, this implies by the work of Auslander,
Ding, and Solberg that $M$ is locally free in codimension one.  We
therefore assume $d \geq 2$.  
The following lemma is
standard. (See \cite[Theorems 16.6, 16.7]{Matsumura}.)

\begin{lemma}\label{finitelength} Let $(R, \m, k)$ be a CM local ring of
dimension at least 2.  Let $X$ be a MCM $R$-module and $L$ a module of
finite length over $R$.  Then $\Ext_R^1(L,X)=0$.\end{lemma}

Recall from Proposition~\ref{triple} that under the assumptions (\ref{exts-d-2}),
the homomorphism $1\otimes \alpha:  \omega \otimes_R M^*\otimes_R M \to
\omega \otimes_R \Hom_R(M,M)$ is injective. 

\begin{lemma}\label{splitmono} If $M$ is locally free on the punctured
spectrum, then the homomorphism $1\otimes \alpha$ is a split
monomorphism with cokernel of finite length.\end{lemma}

\begin{proof} We have the following exact sequence:
\begin{equation}\label{alphaC}
\CD 
0 @>>> \omega \otimes_R M^*\otimes_R M @>1\otimes\alpha>> \omega \otimes_R
\Hom_R(M,M) @>>> C @>>> 0. 
\endCD
\end{equation}
Since $M$ is locally free on the punctured spectrum, $1 \otimes \alpha$
is an isomorphism when localized at any nonmaximal prime of $R$, which
forces $C$ to have finite length.  Since $\omega\otimes_R M^* \otimes_R
M$ is MCM by Proposition~\ref{triple}, $\Ext_R^1(C,\omega\otimes_R M^*
\otimes_R M)=0$, and so (\ref{alphaC}) splits.\end{proof}

\begin{proof}[Proof of Theorem~\ref{reducedmain}]  We will proceed by
induction on $d =\dim R$.  As mentioned above, the case $d=1$ follows
from \cite[Proposition~1.9]{AusDingSolberg}, so we may assume $d \geq 2$, and
that the statement is true for all modules over CM local rings matching
our hypotheses (\ref{exts-d-2}) and having dimension less than that of
$R$.  In particular, we may assume that $M$ is locally free on the
punctured spectrum.  Also, we may assume that $M$ is indecomposable.

First assume that $R$ is Gorenstein. Then $\alpha: M^*\otimes_R M \to
\Hom_R(M,M)$ must be a split monomorphism with cokernel of finite length, by
Lemma~\ref{splitmono}. Since
$\Hom_R(M,M)$ is torsion-free, this implies $\alpha$ is an isomorphism,
and hence that $M$ is free.

Next assume that $R$ is not necessarily Gorenstein, but contains the
rationals.
Consider the following diagram involving the trace homomorphism 
(see Appendix~\ref{trace}).
\begin{equation*}\xy{\xymatrix{
\omega \otimes_R M^* \otimes_R M \ar^-{1\otimes\alpha}[r] \ar_-{1
\otimes \ev}[dr]
& \omega \otimes_R\Hom_R(M,M) \ar^-{1 \otimes \tr}[d]
\\
& \omega\otimes_R R 
}}\endxy\end{equation*}
By Lemma~\ref{Schur}, the diagram commutes.  
Furthermore, by Lemma~\ref{splitmono}, $1 \otimes\alpha$ is a split monomorphism 
with finite-length cokernel $C$, so $\omega \otimes_R \Hom_R(M,M)$ has $C$ as
a direct summand and $1\otimes\alpha$ is surjective onto the complement.  Since 
$\omega$ is torsion-free, $1\otimes \tr$ must kill $C$.  

As $R$ contains $\Q$, $\rank M$ is invertible and so $\tr$ is surjective by 
Corollary~\ref{tronto}.  It follows that the composition $1 \otimes \tr\alpha$ 
is surjective, so that $1\otimes \ev$ is as well.
In other words, the evaluation map $M^*\otimes_R M \to R$ induces a
surjection when tensored with $\omega$.
By Nakayama's Lemma, then, the evaluation map is surjective, and it follows that 
$M$ has a free direct summand.  Since $M$ is indecomposable, $M$ is free.
\end{proof}

\appendix
\section{The Trace of a Module}\label{trace}

In this section we give a general description of the trace of a module. 
Our treatment is intrinsic to the module, and it satisfies the usual
properties of a trace defined for torsion-free modules over a normal domain. 
We include full proofs for convenience.

Throughout this section, let $R$ be a Noetherian ring with total
quotient ring $K$; that is, $K$ is obtained from $R$ by inverting all
nonzerodivisors.  Let $M$ be a torsion-free $R$-module.  The {\it trace
of $M$} will be a certain homomorphism $\tr: \Hom_R(M,M) \to R.$  To
define the trace, let 
$$\alpha: M^* \otimes_R M \to \Hom_R(M,M)$$
be the natural homomorphism defined by $\alpha(f \otimes x)(y) = f(y)
 \cdot $.  Note that dualizing $\alpha$ gives a homomorphism $\alpha^*$ from
$\Hom_R(M,M)^* = \Hom_R(\Hom_R(M,M),R)$ to $(M^* \otimes_R M)^* \cong
\Hom_R(M^*,M^*)$.  It is known (see \cite{Leuschke:gormods}, for
example) that $\alpha$ is an isomorphism if and only if $M$ is free.

\begin{defn} Assume that $\alpha^* : \Hom_R(M,M)^* \to \Hom_R(M^*,M^*)$
is an isomorphism.  The {\em trace of $M$} is defined by $\tr =
(\alpha^*)^{-1}(1_{M^*})$.  We say in this case that $M$ {\em has a
trace.}\end{defn}

Observe that the target of $\alpha^*$ is $(M^* \otimes_R M)^*$, which we
have used $\Hom$-Tensor adjointness to identify with $\Hom_R(M^*,M^*)$. 
Under this identification, the identity map $M^* \to M^*$ corresponds to
the evaluation map $\ev: M^*\otimes_R M \to R$ defined by $\ev(f\otimes
x)= f(x)$. To see this, recall that the $\Hom$-Tensor morphism 
$\Phi_{ABC} : \Hom(A\otimes B, C) \to \Hom(A,\Hom(B,C))$
is defined by $[\Phi_{ABC} (f)(a)](b) = f(a\otimes b)$ for $a \in A$, $b
\in B$.  Taking $A = M^*$, $B=M$, $C=R$, we see that for $x\in M$ and $f
\in M^*$,
$[\Phi_{M^*MR}(\ev)(f)](x) = \ev(f \otimes x) = f(x).$
So $\Phi_{M^*MR}(\ev)$ is the map $M^* \to M^*$ taking $f$ to $f$.  In particular, we could also define the trace by $\tr = (\alpha^*)^{-1}(\ev)$.

Our first proposition generalizes the standard fact that a torsion-free
module over a normal domain has a trace.  

\begin{prop} If $M_\p$ is a free $R_\p$-module for all primes $\p$ of
height one in $R$, and $R$ satisfies Serre's condition (S$_2$), then $M$
has a trace.\end{prop}

\begin{proof} We must show that $\alpha^* : \Hom_R(M,M)^* \to
\Hom_R(M^*,M^*)$ is an isomorphism.  Let $L = \ker(\alpha)$,
$I=\im(\alpha)$, $C= \coker(\alpha)$.  Then dualizing $\alpha$ gives two
exact sequences:
\begin{equation*}\begin{split}
0 \to I^* \to &\Hom_R(M^*,M^*) \to L^* \\
0 \to C^* \to &(\Hom_R(M,M)^* \to I^* \to \Ext_R^1(C,R)
\end{split}
\end{equation*}
Since $\alpha$ is an isomorphism at all minimal primes of $R$, the
annihilator of $L$ is not contained in any minimal prime.  Hence $L$ is
a torsion module, and so $L^* = 0$.  

Since, further, $\alpha$ is an isomorphism at all primes of height one
in $R$, the annihilator of $C$ is not contained in any height-one
prime.  By the assumption that $R$ satisfies condition (S$_2$), then,
$\grade(\Ann C) \geq 2$, so $C^* = \Ext_R^1(C,R) = 0$.  This shows that
$\alpha^*$ is an isomorphism.
\end{proof}

\begin{lemma}\label{freetrace} For $f\in \Hom_R(R^n,R^n)$,  $\tr(f)$ is
the sum of the diagonal entries of a matrix representing $f$.
\end{lemma}

\begin{proof}  Since $R^n$ is free, $\alpha$ is an isomorphism already,
and of course $\alpha^*$ is as well.  Write 
$f = \alpha(\sum_{i=1}^n a_{ij} g_j \otimes e_i),$ where $e_i$ and $g_i$
are the canonical bases for $R^n$ and its dual, respectively.   Then
since $g_j(e_i)=\delta_{ij}$, we see that 
$$\tr(f) = \ev(\sum_{i=1}^n a_{ij} g_j \otimes e_i) = \sum_{1 \leq i,j \leq n} a_{ij}
g_j(e_i) = \sum_{1 \leq j \leq n} a_{jj},$$ 
as desired.
\end{proof}

Recall that the torsion-free $R$-module $M$ is said to have constant
rank $n$ if $K\otimes_R M$ is a free $K$-module of rank $n$.  If this is
the case, we fix a basis $\{e_1, \ldots, e_n\}$ for $K \otimes_R M$, and
let $\{g_1, \ldots, g_n\}$ be the dual basis, so that $g_i(e_j) =
\delta_{ij}$.

\begin{lemma}\label{formulas} Assume that $M$ is a torsion-free
$R$-module of constant rank and that $M$ has a trace.  Then for any $f
\in M^*$ and $x \in M$, we have 
$x = \sum_{i=1}^n g_i(x) e_i$
and
$\tr(f) = \sum_{i=1}^n g_i(\hat{f}(e_i)),$
where $\hat{f} = K \otimes_R f$.
\end{lemma}

\begin{proof} Since $M$ is torsion-free, it embeds into a free
$R$-module and so the homomorphism $M \to K\otimes_R M$ is injective. 
Considering $x$ as an element of $K \otimes_R M$, write 
$x = \sum_{j=1}^n a_j e_j,$
where the $a_j$ are elements of $K$.  Then 
a short computation using the definition of the $g_i$ shows that
$\sum_{i=1}^n g_i(x) e_i =x$.  For the other assertion, pass to the
total quotient ring $K$.  Since $K\otimes_R M$ is free, Lemma~\ref{freetrace} implies that the trace of $\hat{f}$ is the sum of the diagonal elements of a matrix $(a_{ij})$ representing $\hat{f}$.  Since
$g_i(\hat{f}(e_i)) = a_{ii}$,  the statement follows.\end{proof}

\begin{cor}\label{tronto} Assume that $M$ is a torsion-free module of constant rank
and has a trace.  If $\rank(M)$ is invertible in $R$, then $\tr$ is
surjective from $\Hom_R(M,M)$ to $R$.\end{cor}

\begin{lemma} \label{Schur} Assume that $M$ is a
torsion-free of constant rank and that $M$ has a trace.  Then we have
$\tr \alpha = \ev$ as homomorphisms from $M^*\otimes_R M$ to
$R$.\end{lemma}

\begin{proof} For any $f \in M^*$ and $x \in M$, a straightforward computation using
Lemma~\ref{formulas} shows that $f(x) = \tr(\alpha(f \otimes x))$.
\end{proof}

\providecommand{\bysame}{\leavevmode\hbox to3em{\hrulefill}\thinspace}
\providecommand{\MR}{\relax\ifhmode\unskip\space\fi MR }
\providecommand{\MRhref}[2]{%
  \href{http://www.ams.org/mathscinet-getitem?mr=#1}{#2}
}
\providecommand{\href}[2]{#2}

\end{document}